\documentclass[12pt]{amsart}
\usepackage{amsthm, amscd}
\usepackage{amssymb,latexsym}

\theoremstyle{plain}
\newtheorem*{introtheorem}{Theorem}
\newtheorem*{introcorollary}{Corollary}
\newtheorem{theorem}{Theorem}[section]
\newtheorem{proposition}[theorem]{Proposition}
\newtheorem{corollary}[theorem]{Corollary}
\newtheorem{lemma}[theorem]{Lemma}

\theoremstyle{definition}

\newtheorem{remark}[theorem]{Remark}
\newtheorem{example}[theorem]{Example}

\theoremstyle{remark}
\newtheorem*{nctheorem}{Theorem}

\newcommand{\lra}{\longrightarrow}
\newcommand{\noi}{\noindent}
\newcommand{\PP}{\mathbf{P}}

\newcommand{\CC}{\mathbf{C}}

\newcommand{\OO}{\mathcal{O}}
\newcommand{\II}{\mathcal{I}}
\newcommand{\JJ}{\mathcal{J}}

\newcommand{\blowup}{\text{Bl}}
\newcommand{\sbl}{\vskip 3pt}

\begin{document}

\title {A Geometric  Effective Nullstellensatz}

\author{Lawrence Ein} 
\address{Department of Mathematics \\
University of Illinois at Chicago \\ 851 South Morgan
St., M/C. 249
\\ Chicago, IL  60607-7045}
\thanks{Research of first author  partially
supported by NSF Grant DMS 96-22540}
\email{ein@math.uic.edu}
\author{Robert Lazarsfeld} 
\address{Department of Mathematics
\\ University of Michigan \\ Ann Arbor, MI  48109}
\email{rlaz@math.lsa.umich.edu}
\thanks{Research of  second author  partially
supported by the J. S. Guggenheim Foundation and NSF
Grant DMS 97-13149 }

\maketitle
\setcounter{section}{-1}


\section*{Introduction}
The purpose of this paper is to present a
geometric theorem which clarifies and extends in
several directions work of Brownawell, Koll\'ar 
and others on the effective Nullstellensatz.
Specifically, we work on an arbitrary smooth complex
projective variety X, with the previous ``classical"
results corresponding to the case when X is projective
space. In this setting we prove a local effective
Nullstellensatz for ideal sheaves, and a corresponding
global division theorem for adjoint-type bundles. We
also make explicit the connection with the intersection
theory of Fulton and MacPherson. Finally, constructions
involving products of prime ideals that appear in
earlier work are replaced by geometrically more natural
conditions involving orders of vanishing along
subvarieties. 

Much of the previous activity in this area has
been algebraic in nature, and seems perhaps not
well-known in detail among geometers. Therefore we have
felt it worthwhile to include  here a rather extended
Introduction. We start with an overview of the
questions and earlier work on them. Then we present the
set-up and statement of our main theorem. We conclude 
with a series of examples (which can be read before the
general result) of what it yields in special cases.

\subsection*{Background}

 In recent
years  there has been a great deal of interest in the
problem of finding effective versions of Hilbert's
Nullstellensatz. The classical theorem of course states
that given polynomials
\[ f_1 , \dots , f_m \in \CC[t_1, \dots, t_n ] , \]
if the $f_j$ have no common zeroes in $\CC^n$  then
they generate the unit ideal, i.e. there exist $g_j \in
\CC[ t_1, \dots , t_n ]$ such that 
\begin{equation}
\sum g_j f_j = 1. \tag{*}
\end{equation}
A first formulation of the problem is to bound the
degrees of the $g_j$ in terms of those of the $f_j$. 
Current work in this area started with a theorem of
Brownawell \cite{Brownawell1}, who showed that if  $\deg
f_j
\le d$ for all $j$, then one can find $g_j$ as in (*)
such that 
\begin{equation}
\deg g_j  \le n^2 d^n + nd. \tag{B1}
\end{equation}
Brownawell's argument was arithmetic and analytic in
nature, drawing on height inequalities from
transcendence theory and the classical theorem of
Skoda. Shortly thereafter, Koll\'ar \cite{Kollar1}
gave a more elementary and entirely algebraic proof of
the optimal statement that in the situation above, one
can in fact take
\begin{equation} \deg ( g_j f_j ) \le d^n \tag{K1}
\end{equation}
 provided that $d \ne 2$.\footnote{Here
and below we are oversimplyfing slightly Koll\'ar's
results. He actually establishes a  more precise
statement allowing for the
$f_j$ to have different degrees, and giving stronger
estimates when $m \le n$. Furthermore, he works over an
arbitrary ground field.} 

Koll\'ar deduces (K1) as an immediate consequence of a
rather surprising theorem in the projective setting.
Specifically, consider a homogeneous ideal $J \subseteq
\CC[\ T_0, \dots, T_n \ ]$. Then of course $J$ contains
 some power of its radical. The main theorem of
\cite{Kollar1} is the effective statement that if $J$ is
generated by forms of degree
$\le d \ (d \ne 2)$, then  already 
\begin{equation}
\big( \sqrt{J} \big ) ^{d^n} \subseteq J. \ \ 
\tag{K2} \end{equation} [Proof of (K1):  let $F_j \in
\CC[T_0, \dots, T_n]$ be the homogenization of
$f_j$. Then the common zeroes of the $F_j$ lie in the
hyperplane at infinity $\{ T_0 = 0 \}$, and
consequently $T_0 \in \sqrt{ (F_1, \dots, F_m )}$.
Therefore $(T_0)^{d^n} = \sum G_j F_j$ thanks to (K2),
and (K1) follows upon dehomogenizing.] By
analyzing Koll\'ar's proof, Brownawell
\cite{Brownawell2} subsequently shed a somewhat more
geometric light on this result. Namely, still assuming
that $J$ is generated by forms of degree $\le d $,
he shows that there exist reduced and
irreducible subvarieties $W_i \subset
\PP^n$ with
\[
 \cup W_i = Z =_{\text{def}} \text{Zeroes}(\sqrt{J}),
\]
plus positive integers $s_i > 0$, satisfying the
following properties. First, one has the degree bound
\begin{equation}
\sum s_i \deg W_i \le d^n, \tag{B2}
\end{equation}
so that in particular $\sum s_i \le d^n$. Secondly, if
$I_{W_i}$ denotes the homogeneous ideal of
$W_i$, then 
\begin{equation} \prod \big ( I_{W_i} \big )^{s_i}
\subseteq J. 
\tag{B3}
\end{equation}
This formulation is referred to as the  ``prime-power
Nullstellensatz" or the ``algebraic Bezout
theorem".\footnote{As explained in 
\cite{Brownawell2} one should take here $W_0 =
\emptyset$, with
$I_{W_0}  = (T_0, \dots , T_n)$, and assign to  $W_0$
``honorary degree" one.} Since $\sqrt{J} \subseteq
I_{W_i}$ for every $i$, it is immediate that (B3) and
(B2) imply (K2), and in fact (B3) improves (K2) unless
every component of   $Z$ is a linear space. However
Brownawell's construction does not provide a clearly
canonical choice for the
$W_i$. We refer to
\cite{Tessier} and \cite{BernsteinStruppa} for 
excellent surveys of this body of work, and to
\cite{Bernstein} for a discussion of some analytic
approaches to these questions. Recently Sombra 
\cite{Sombra} proved an analogue of (K2) for
projectively Cohen-Macaulay varieties $X \subset
\PP^N$, from which he deduces an interesting
generalization of (B1) for sparse systems of
polynmials (see Examples 2 and 3 below). Motivated in
part by Sombra's work, Koll\'ar
\cite{Kollar2} has generalized  these results
to arbitrary ideals in the polynomial ring.

While this picture is fairly complete from an
algebraic point of view, a number of geometric
questions present themselves. First, it is natural to
ask whether the results of Koll\'ar and Brownawell  ---
which involve homogeneous ideals in the polynomial ring
--- can be seen as the case $X = \PP^n$ of a more
general picture involving an arbitrary smooth
projective variety $X$:  Sombra's theorem gives one
step in this direction.  Next, one might hope to clarify
the connection with intersection theory that is
evidently lurking here. Finally, it is
difficult geometrically to determine whether a given
polynomial lies in a product of ideals, and from this
point of view one would like to replace the product
of prime powers occuring in (B3) by an
intersection of symbolic powers  defined by orders of
vanishing along subvarieties. The theorem we present in
this paper attemps to address these
questions.\footnote{We should state at the outset
however that in the ``classical" case $X = \PP^n$ our
numerical bounds are in some instances slightly weaker
than those of Koll\'ar-Brownawell.}

\subsection*{Set-up and Statements.} Turning to a 
detailed presentation of our results, we start by
introducing the set-up in which we shall work, and by
fixing some notation.  Let
$X$ be a smooth complex projective variety of dimension
$n$, and let 
\[
D_1, \dots, D_m \in |D|
\]
be effective divisors on $X$ lying in a given linear
series. Set $L = \OO_X(D)$, and let $s_j \in \Gamma(X,
L)$ be the section defining $D_j$. We denote by $B$ the
scheme-theoretic intersection
\[
B = D_1 \cap \dots \cap D_m \subset X, 
\]
and we let 
\[
\JJ = \sum \OO_X(-D_j) \subset \OO_X  
\]
be its ideal sheaf. Finally,  set $Z =
B_{\text{red}}$, so that $Z = \text{Zeroes}(\sqrt{\JJ})$
is the reduced scheme defined by the radical of $\JJ$.

Recall next from \cite{Fulton}, Chapter 6, \S 1, that
the scheme $B$ canonically determines a
decomposition
\[ Z = Z_1 \cup \dots \cup Z_t \]
of $Z$ into (reduced and irreducible) {\it distinguished
 subvarieties} $Z_i \subset Z$,
together with positive integers  $r_i > 0$. We will
review the precise definition in \S 2, but for the
moment suffice it to say that the
$Z_i $ are the supports of the irreducible
components of the projectivized normal cone 
$\PP(C_{B/X})$  of $B$ in
$X$. The coefficient $r_i$ attached to
$Z_i$ arises as the multiplicity of the corresponding
component of the exceptional divisor in the
(normalized) blowing up of $X$ along $B$. Every
irreducible component of $Z$ is distinguished,
but there can be  ``embedded" distinguished subvarieties
as well.  We denote
by
$\II_{Z_i} \subseteq \OO_X$ the ideal sheaf of $Z_i$,
and by
$\II_{Z_i}^{<r>}$ its
$r^{\text{th}}$ symbolic power, consisting of germs of
functions that have multiplicity
$ \ge r$ at a general point of $Z_i$.
\sbl
Our main result is the
following:
\begin{introtheorem}
With notation and assumptions as above, suppose that $L$
is ample.
\begin{enumerate} 
\item[(i).] The distinguished subvarieties $Z_i
\subset X$ satisfy the degree bound
\[
 \sum r_i \cdot \deg_L(Z_i) \le deg_L(X) = \int_X
c_1(L)^n, 
\]
where as usual the $L$-degree of a subvariety $W
\subseteq X$ is the integer $\deg_L(W) = \int_W
c_1(L)^{\dim(W)}$.

\sbl

\item[(ii).] One has the inclusion
\[
\II_{Z_1}^{<n \cdot r_1>} \cap \dots \cap
\II_{Z_t}^{<n \cdot r_t>}
\subseteq \JJ. \]
In other words, in order that a function (germ) $\phi$
lie in $\JJ$, it suffices that $\phi$ vanishes to order
$\ge n r_i$ at a general point of each of the
distinguished subvarieties $Z_i$.
\sbl
\item[(iii).] Denote by $K_X$ a canonical divisor
of $X$, and let
$A$ be  a divisor on
$X$ such that
$A -(n+1) D$ is ample. If 
\[s \in
\Gamma(X, \OO_X(K_X + A ))\]
 is a section which
vanishes to order $\ge (n+1) \cdot  r_i$ at the general
point of each
$Z_i$, then one can write
\[
s = \sum s_j h_j \quad \text{for some sections} 
\ h_j \in \Gamma(X,
\OO_X(K_X + A - D_j)),
\]
where as above $s_j \in \Gamma(X, \OO_X(D_j))$ is the
section defining $D_j$.
\end{enumerate}
\end{introtheorem}

 As in
Brownawell's algebraic Bezout theorem, the
inequality in (i)  serves in effect to bound the
coefficients $r_i$ from above. One should view (ii) as
a local effective Nullstellensatz. Together with (i) it
immediately implies  the first statement of
the
\begin{introcorollary} 
\begin{enumerate}
\item[(a).] With notation and assumptions as
above:
\[ \big( \sqrt{\JJ} \big) ^{n
\cdot \deg_L(X)}
\subseteq \big(\sqrt{\JJ}
\big) ^{n \cdot \max\{r_i\}} \subseteq \JJ. 
\]
More generally,
\[
\big( \sqrt{ \JJ} \big) ^{< n \cdot \deg_L(X) >}
\subseteq \JJ, 
\]
where the symbolic power on the left  denotes
the sheaf of all functions that vanish to at least  the
indicated order at every point of $Z$. 
\sbl
\item[(b).] If $s \in \Gamma(X,\OO_X(kD))$ is a section
which has multiplicity $\ge (n+1) \int c_1(L)^n$ at
every point of $Z$, then if $k \gg 0 $ is sufficiently
large there exist   $h_j \in \Gamma\big( X,
\OO_X((k-1)D) \big)$ such that $s = \sum s_j
h_j$. \qed
\end{enumerate}
\end{introcorollary}
\noi It is perhaps already somewhat surprising  that
there are tests for membership in an ideal that depend
only on orders of vanishing along its zero-locus. Note
that the Theorem applies to an arbitrary ideal sheaf
$\JJ$ as soon as $L
\otimes
\JJ$ is globally generated. So from a qualitative point
of  view one may think of the Corollary as giving global
constraints on the local complexity of $\JJ$. On the
quantitative side, we remark that the factor of $n$
appearing in (ii) and statement (a) of the Corollary
can be replaced by
$\min(m,n)$, and similarly in (iii) and (b) one can
substitute $\min(m, n+1)$ for
$(n+1)$. The results of Koll\'ar and Brownawell  might
suggest the hope that one could drop  these factors
altogether, but examples (see \ref{counterexample})
show that this is not possible, at least with the $Z_i$
and
$r_i$ as we have defined them.  However it
is possible that  (a) holds with the
exponent $n \cdot \deg_L(X)$ replaced by $\deg_L(X)$,
with an analogous improvement of (b).

The proof of the Theorem is quite elementary and, we
hope, transparent. It consists of three steps. First
(\S 1) we use vanishing theorems to give a simple
algebro-geometric proof of a statement of Skoda type.
The theorem in question establishes local and global
criteria involving  some multiplier-type ideal sheaves
$\II_\ell$ to guarantee that one can write a given germ
$\phi \in \OO_X$  or global section $s \in \Gamma(X,
\OO_X(K_X + A))$ in terms of the $s_j \in \Gamma(X, L)$.
The local statement was originally proved in
\cite{Skoda} using  $L^2$-methods, and while Skoda's
result is  well known in analytic geometry and
commutative algebra (cf. \cite{LT},
\cite{Huneke} and \cite{Lipman}), it seems to be less
familiar to algebraic geometers. We hope
therefore that the discussion in \S 1 -- which in
addition contains an extension of these results to
higher powers of
$\JJ$ -- may be of independent interest.\footnote{In
fact   one can deduce the local effective
Nullstellensatz  directly from the theorem of Brian\c
con and Skoda for regular local rings (Remark
\ref{rtt}). From our perspective however the
local and global statements are two sides of the same
coin, and in essence we end up reproving  Brian\c
con-Skoda.} The next point
(\S 2)  is to relate the  sheaves $\II_\ell$ to orders
of vanishing along the Fulton-MacPherson distinguished
subvarieties
$Z_i$. Section 2 also contains a geometric
characterization of these subvarieties, in the spirit
of van Gastel, Flenner and Vogel
(\cite{VanGastel},\cite{FlennerVogel}). Finally,  a
simple calculation of intersection numbers gives the
degree bound (\S3).  It is interesting to observe that
while  the final outcome is quite different,
essentially all of these techniques have antecendents
in earlier work in this area.

\subsection*{Examples} Finally, in order to give a
feeling for the sort of concrete statements that come
out of the Theorem, we  conclude this Introduction
with a few examples. 
 
\noi \textbf{Example 1. }
Consider the ``classical" case $X =
\PP^n$ and
$L = \OO_{\PP^n}(d)$, so that we are dealing with $m$
homogeneous polynomials 
\[
s_1, \dots , s_m \in \CC[T_0, \dots T_n]
\] 
of degree $d$. Then the degree bound in part (i) of the
Theorem says that
\begin{equation}
\sum r_i \cdot d^{\dim(Z_i)} \cdot \deg(Z_i) \le d^n, 
\tag{*} \end{equation}
where here $\deg(Z_i)$ is the standard degree (with
respect to $\OO_{\PP^n}(1)$). 
The conclusion of statement (iii)  is that
if $s$ is a homogeneous polynomial of degree
$\ge (n+1)(d-1) + 1$ vanishing to order $\ge r_i(n+1)$
on each of the $Z_i$ then $s$ lies in the homogeneous
ideal $J$ spanned by the $s_j$. In other words, if
$I_{Z_i}^{<r>}$ denotes the homogeneous primary
 ideal of all polynomials
having multiplicity $\ge r$ at a general
point of $Z_i$, and if $(T_0, \dots , T_n)$
denotes the irrelevant maximal ideal, then we have
\[
(T_0, \dots, T_n)^{(dn + d - n)} \cap
I_{Z_1}^{<(n+1)r_1>}
\cap
\dots
\cap  I_{Z_t}^{<(n+1)r_t>} \subset J.
\]
By analogy with Brownawell's ``prime-power"
formulation of Koll\'ar's theorem,
one might think of this as a ``primary decomposition"
version of the Nullstellensatz. Comparing this with
Brownawell's  statement (B3), the most
surprising difference is that one can ignore here any of
the ``embedded" distinguished subvarieties $Z_i$ for
which the corresponding coefficient $r_i$ is small.
Numerically, the factor of 
$d^{\dim(Z_i)}$ in (*) strengthens   (B2), but the
factor of $(n+1)$ in the exponent prevents one from
recovering (K2) in the  ``worst" cases when every
component of $Z$ has small degree.

\noi \textbf{Example 2.}   M. Rojas has observed that
following the model of  \cite{Sombra} one can apply the
Theorem to suitable toric compactifications $X$ of
$\CC^n$ to obtain extensions of the results (B1)
and (K1) of Brownawell and Koll\'ar to certain sparse
systems of polynomials (see also \cite{KK}, \cite{GKZ},
and
\cite{Rojas1} for other other applications of toric
geometry to sparse systems of polynomials). In some
settings, the numerical bounds that come out strengthen
Sombra's.  We refer to the forthcoming preprint
\cite{Rojas2} of Rojas for the precise  statements, but
illustrate their flavor in a special case. Consider as
above polynomials  $f_j \in \CC[t_1, \dots, t_n]$ and
suppose that one is given separate degree bounds in
each of the variables $t_k$:
\[ \deg _{t_k}(f_j) \le d_k \ \ \forall \  j. \]
Assuming that the $f_j$ have no common zeroes in
$\CC^n$,  then one can find $g_j$ with $ \sum g_j
f_j = 1$ where now
\begin{equation} \deg (g_j f_j )  \ \le (n+1) !  \ \ d_1
\cdot \dots \cdot d_n. \tag{*} \end{equation} 
(By way of comparison, Sombra's general theorem
yields in this setting the analogous inequality with
the factor of $(n+1)!$ replaced by $n^{n+3}$.)  If for
instance one thinks of
$d_1,
\dots , d_{n-1}$ as  being fixed, then (*) gives a
linear bound in the remaining input degree $d_n$.
[To prove (*), one applies the Theorem to $X = \PP^1
\times \dots \times \PP^1$ and $L = \OO(d_1, \dots,
d_n)$, and argues as in the proof that (K2) implies
(K1).]

\noi \textbf{Example 3.} Our last example is a variant
of a result of Sombra \cite{Sombra}, (1.8). In
the situation of the Theorem, suppose that
$H$ is a very ample divisor on $X$ which is sufficiently
positive so that that $H - K_X - (n+1)D$ is ample, and
consider the embedding
$X \subset \PP^N = \PP$ defined by the complete linear
system $|H|$. 
Let $I \subset S := \CC[T_0, \dots, T_N]$ be
the homogeneous ideal of $X$ under this embedding,
let $R = S/I$ be the homogeneous coordinate ring of
$X$, and let
$F_1, \dots, F_m \in R$ be homogeneous elements of
degrees $\le d$. Let $P \in R$ be a homogeneous element
lying in the radical of the ideal $(F_1, \dots, F_m)$.
Then 
\[ P^{ (n+1)d^n \deg X} \in (F_1, \dots, F_m), \]
where $\deg X = (H^n)$ denotes the degree of $X$ in the
projective embedding defined by $|H|$. (As in
\cite{Sombra}, one first reduces to the case where all
the $F_j$ are of equal degree $d$.) When $I$ is a
Cohen-Macaulay ideal --- which for sufficiently
positive $H$ is equivalent to the vanishings $H^i(X,
\OO_X) = 0$ for $0 < i < \dim X$ --- this is a slight
numerical improvement of Sombra's
result (which however does not require the variety
defined by $I$ to be non-singular).

\sbl
We wish to thank W. Fulton, J. Koll\'ar, M. Rojas,  K.
Smith and M. Sombra for valuable correspondence and
discussions.

\sbl
\section{Notation and Conventions}
\noi {(0.1).} We work throughout with varieties
and schemes defined over the complex numbers.

\noi {(0.2).} Let $X$ be a smooth variety, and
$\phi \in \OO_X$ the germ of a regular function defined
in the neighborhood of a point $x \in X$. We say that
$\phi$ {\it vanishes to order $\ge r$ at $x$}, or that
$\phi$ {\it has multiplicity $\ge r$ at $x$} if $\phi
\in m_x^r$, where $m_x \subset \OO_xX$ is the maximal
ideal of $x$. Equivalently, all the partials of $\phi$
of order $< r$ should vanish at $x$.  If $Z \subset X$
is an irreducible subvariety, with ideal sheaf $\II_Z
\subset \OO_X$, we denote by $\II_Z^{<r>} \subset \OO_X$
the sheaf of germs of functions that vanish to order
$\ge r$ at a general (and hence at every) point of $Z$.
It is a theorem of Nagata and Zariski
(cf. \cite{Eisenbud}, Chapter 3, Section 9) that this
coincides with the
$r^{\text{th}}$ symbolic power of $\II_Z$, although
there is no loss here in taking this as the definition
of symbolic powers. Evidently $\II_Z^r \subset
\II_Z^{<r>}$, but when $Z$ is singular the inclusion
may well be strict. 

\noi {(0.3)} Let $X$ be a smooth projective
variety of dimension $n$. A line bundle $L$  on
$X$ is {\it numerically effective} (or {\it nef}) if 
\[
\int_C c_1(L) \ge 0
\]
 for every irreducible curve $C \subset X$.
A fundamental theorem of Kleiman (cf.
\cite{Hartshorne}, Chapter I, \S 6) implies that any
intersection number  involving the product of Chern
classes of nef line bundles with an effective cycle is
non-negative.  A nef line bundle is {\it big} if its top
self-intersection is strictly positive:
\[ \int_X c_1(L)^n > 0. \]  For a divisor $D$ on $X$,
we define nefness or bigness by passing to the
associated line bundle $\OO_X(D)$.

\noi {(0.4)} The basic global vanishing theorem we
will use is the following extension by Kawamata and
Viehweg of the classical Kodaira vanishing theorem:
\begin{nctheorem} Let $X$ be a smooth complex
projective variety, and let $K_X$ denote a canonical
divisor on $X$. If $D$ is a big and nef divisor on $X$
then
\[ H^i(X, \OO_X(K_X + D)) = 0 \ \ \text{for } \ i>0.
\]
\end{nctheorem}
\noi One of the benefits of allowing merely big and nef
bundles is that this result then implies a local
vanishing theorem for higher direct images. For our
purposes, the following statement will be sufficient:
\begin{nctheorem} Let $X$ be a smooth
quasi-projective complex variety, and let $f : X \lra
Y$ be a generically finite and surjective projective
morphism. Suppose that $D$ is a divisor on $X$
which is nef for $f$, i.e. whose restriction to every
fibre of $f$ is nef. Then
\[R^j f_* \big ( \OO_X(K_X + D) \big )  = 0 \ \
\text{for } \ j >0. \]
\end{nctheorem}
\noi This  is called vanishing
for the map $f$. We refer to \cite{Kollar0} for a very
readable introduction to the circle of ideas
surrounding vanishing theorems, and to \cite{KMM},
(0.1) and (1.2.3) for a more technical and detailed
discussion. 
\section{A Theorem of Skoda Type}

In this section we use vanishing for big and nef line
bundles to give a simple algebro-geometic proof of a
theorem of Skoda type. In his  classical paper
\cite{Skoda}, Skoda  uses $L^2$ techniques to establish
an analytic criterion guaranteeing that a
germ
$f \in \CC\{z_1, \dots , z_n\}$ lies in the ideal
generated by a given collection of functions
$f_1, \dots, f_m \in \CC\{z_1, \dots , z_n\}$. In view
of the close connection that has emerged in recent years
between such
$L^2$ methods and vanishing theorems (cf.
\cite{Demailly} for a survey), it is natural to
expect that one can recover statements of
this sort via vanishing. We carry this out here.
Besides being very elementary and transparent, the
present approach has the  advantage of
simultaneously giving global results. A special case
of Skoda's theorem played a role in Siu's recent work
\cite{Siu} on the deformation invariance of
plurigenera, and it was algebrized as below by
Kawamata \cite{Kawamata}. 

 Let $X$ be a smooth irreducible quasi-projective
complex  variety of dimension $n$. We emphasize that
for the time being
$X$ need not be projective, and  in fact for the local
results one might want to think of $X$ as
representing    the germ of an algebraic variety. Let
\[
\JJ \subseteq \OO_X
\]
be an ideal sheaf defining a proper subscheme $B \subset
X$.  For each $\ell \ge 1$ we associate to $\JJ$ a
multiplier-type ideal sheaf 
\[ \II_\ell \subset \OO_X \]
as follows. Start by forming the blow-up
\[ \nu_o : V_0 = \blowup_B(X) \lra X
\]
of $X$ along $B$, and then take a resolution of
singularities $Y \lra V_0$ to get a birational map
\[ f : Y \lra X. \]
The ideal $\JJ$ becomes principal on $V_0$ and hence
also on $Y$. More precisely, let $F$  be the pull-back
to $Y$ of the exceptional divisor on $V_0$.  Then
\[
\JJ \cdot \OO_Y = \OO_Y(-F).
\]
We set 
\[
\II_{\ell} = f_* \big ( \OO_Y( K_{Y/X} - \ell F ) \big)
,
\] 
where $K_{Y/X} = K_Y - f^* K_X$ is the relative
canonical divisor of $Y$ over $X$. Note that
$\II_\ell \subseteq f_*  \OO_Y( K_{Y/X}  )
 = \OO_X$, so that $\II_\ell$ is indeed an ideal sheaf
on $X$. One can
check by standard arguments that it is independent of
the choice of a resolution, although we don't actually
need this fact here. In the setting of local algebra,
such ideals were introduced and studied by Lipman
\cite{Lipman}. One could also define
$\II_\ell$ via an $L^2$ integrability condition, as in
Skoda's  paper \cite{Skoda}. We refer to \cite{Ein}
for a discussion, from an algebro-geometric viewpoint,
of multiplier ideals of this sort.

Our object is to relate the ideal sheaves $\II_\ell$ to
$\JJ$. To this end let
$L$ be any line bundle on
$X$ such that $L \otimes \JJ$ is globally generated.
Choose global sections
\[ 
s_1, \dots , s_m \in \Gamma(X, L \otimes \JJ)
\]
generating $L \otimes \JJ$, and set 
$
D_j = \text{div}(s_j) \in |L|.
$
Thus the subscheme $B \subset X$ defined by $\JJ$ is
just the scheme-theoretic intersection of the $D_j$.
Note that all the $D_j$ are linearly equivalent: for
convenience we will sometimes write $D$ for any divisor
in their linear equivalence class. Since $L \otimes
\JJ$ is globally generated, so is its inverse image
\[ N =_{\text{def}} f^*L \otimes \OO_Y(-F).
\]
In fact, we can write
\[ f^* D_j = F + D_j^\prime \]
where the $D_j^\prime \in |N|$ are effective divisors
on $Y$ that generate a base-point free linear system.

Pushing forward the evident map
\[
\OO_Y(K_{Y/X} - (\ell -1) F) \otimes f^*
L^* = \OO_{Y}(K_{Y/X} - \ell F - D_j^{\prime})
\overset{\cdot D_j^{\prime}}{\lra} \OO_{Y}(K_{Y/X} -
\ell F ),
\] 
determines a sheaf homomorphism 
\[ \sigma_j : \II_{\ell - 1} \otimes L^* \lra \II_\ell
\] 
on $X$. 
Observe that $\sigma_j$ is induced by multiplication by
$s_j$ in the sense that one has a commutative diagram
\[
\begin{CD}
\II_{\ell -1} \otimes L^* @>\sigma_j>> \II_\ell \\
@VVV @VVV \\
L^* @>{\cdot s_j}>> \OO_X
\end{CD},
\]
where the vertical maps arise from the natural
inclusions of $\II_{\ell -1}$ and $\II_{\ell}$ in 
$\OO_X$. This may be verified by pushing forward the
corresponding commutative square
\[
\begin{CD}
\OO_Y(K_{Y/X} - (\ell -1) F) \otimes f^*
L^* @>\cdot D_j^\prime>> \OO_{Y}(K_{Y/X} -
\ell F ) \\
@V{\cdot (\ell -1) F}VV @VV{\cdot \ell F}V \\
\OO_Y(K_{Y/X}) \otimes f^* L^* @>{\cdot f^* D_j}>>
\OO_Y(K_{Y/X})
\end{CD}
\]
 of invertible sheaves
on $Y$. In particular, the image of $\sigma_j$ lies in
the ideal sheaf $\OO_X(-D_j)$ of $D_j$.

We now come to the main result of this section:
\begin{proposition} \label{T1} 
\begin{enumerate}
\item[(i).] \textnormal{(Skoda's Theorem, cf.
\cite{Skoda}, \cite{Lipman}.)} If $\ell \ge \min(m,n)$
then the sheaf homomorphism
\[ 
\sigma =_{\text{def}} \sum_{j = 1}^m \sigma_j :
\bigoplus_{j = 1}^{m}
\
 \II_{\ell -1} \otimes L^* {\lra}
\II_\ell
\]
is surjective. In particular, \[\II_\ell \subset \JJ.\]

\item[(ii).] Assume that $X$ is
projective, and fix  $\ell  \ge \min(m , n+1)$.
Let $A$ be a divisor on
$X$ such that
$A  - \ell D$ is ample (or big and nef). Then the  map
on global sections
\[ 
\bigoplus_{j =1}^m H^0(X, \OO_X(K_X + A - D_j) \otimes
\II_{\ell -1} ) \lra H^0(X, \OO_X(K_X + A) \otimes
\II_\ell) \]
induced by $\sigma$ is surjective. In particular if 
\[ s
\in H^0(X,
\OO_X(K_X + A))\]
 lies in the subspace
$
H^0(X, \OO_X(K_X + A)\otimes \II_{\ell}) \subseteq
H^0(X, \OO_X(K_X + A))
$, then
\[ s = \sum h_j s_j \qquad \text{for some} \ 
h_j \in H^0(X, \OO_X(K_X + A - D_j)).
\]
\end{enumerate}

\end{proposition} 

\begin{proof} As in \cite{LT}, \S 5, we argue via a
Koszul complex. Working on
$Y$, let
$P$ be the  vector bundle
\[ 
P = \bigoplus_{j=1}^m \OO_Y(-D_j^\prime) \cong
\bigoplus_{j=1}^m N^*. 
\]
Then the $D_j^\prime$ determine in the evident way
a surjective homomorphism
$P \lra \OO_Y$. Form the corresponding Koszul complex
and for fixed $\ell$ twist by $Q = Q_\ell =_{\text{def}}
\OO_Y(K_{Y/X} - \ell F)$:
\begin{equation}
\dots \lra \Lambda^2 P \otimes Q \lra P \otimes Q \lra
Q \lra 0. \tag{*} \end{equation}
For (i) we need to establish the surjectivity of the
push-forward homomorphism:
\[ 
\begin{CD}
  f_* \big ( P \otimes Q
\big ) @>>> f_* Q \\
@| @| \\
\bigoplus \II_{\ell-1} \otimes L^* @. \II_\ell
\end{CD}.
\]
Chasing through the  exact sequence (*), we see
that it is enough to establish the vanishings:
\begin{equation}
R^j f_* \big ( \Lambda^{j+1} P \otimes Q \big )= 0
\ \ \text{for} \ 1 \le j \le n.
\tag{**} \end{equation}
Since all the fibres of $f$ have dimension $\le n - 1$,
 the vanishing of the $n^{\text{th}}$ direct image $R^n
f_*$ in (**) is free. So  we can limit attention
to $ j
\le n-1$. Furthermore, as $P$ has rank $m$,
(**) is trivial if $j + 1 > m$. Thus all told we are
reduced to considering only $ j+1 \le \min(m,n)$ in
(**). 

 Now
\[
\Lambda^i P \otimes Q = \Lambda^i P \otimes
\OO_Y(K_{Y/X} - \ell F) \cong
\oplus \OO_Y(K_{Y} \otimes N^{\otimes (\ell - i)} )
\otimes  f^* \OO_X(-K_X) \otimes f^* L^{\otimes -\ell}.
\]
But $N$ is globally generated, and hence is nef for $f$
(and globally nef when $X$ is projective). Furthermore, 
thanks to the projection formula twisting by bundles
pulling back from
$X$ commutes with taking higher direct images. Hence it
follows from vanishing for
$f$ (cf. (0.4)) that one has the vanishing of all
the higher direct images
\[
R^j f_* \big ( \Lambda^i P \otimes Q \big )= 0
\ \ \text{for} \ j > 0, \ i \le \ell.
\]
This proves (**) (when  $j+1
\le \min(m,n)$ and $\ell \ge \min(m,n)$), and with
it statement (i).
 
The second assertion
follows similarly by applying global vanishing for big
and nef divisors on $Y$. In fact, twisting  by $f^*
\OO_X(K_X + A)$, we need to prove the surjectivity of
the homomorphism
\[ 
\begin{CD}
  H^0 \big( Y,  P \otimes Q \otimes f^* \OO_X(K_X
+ A)
\big ) @>>> H^0 \big ( Y,  Q \otimes f^* \OO_X(K_X +
A) \big )
\\ @| @| \\
\bigoplus H^0 \big( X, \OO_X(K_X + A - D_j) \otimes
\II_{\ell -1} \big) @.
H^0\big( X, \OO_X(K_X + A) \otimes \II_\ell \big )
\end{CD}.
\]
determined by the map on the right in (*). Chasing
again through that sequence it suffices to establish the
vanishings
\begin{equation} 
H^j \big( Y, \OO_Y(K_Y + f^*(A -\ell D)) \otimes
N^{\otimes (\ell - j -1)} \big) = 0 \ \text{for} \ 0
< j \le
 \min(m-1, n). \tag{***}
\end{equation}
But by hypothesis $f^*(A - \ell D)$ is big and nef, and
$N$ is nef. So provided  that $\ell \ge \min(m, n+1)$
the bundle occuring in (***) is big and nef, and we are
done thanks to (0.4). 
\end{proof}

\vskip 5pt
\smaller[1]

Although not required for the main development, as in 
\cite{LT}, \cite{Lipman} and \cite{Huneke}, Chapter 5,
it is of some interest to extend these results to higher
powers of $\JJ$. To this end, given a multi-index $J =
(j_1,
\dots , j_m)$ of length $|J| = \sum j_\alpha = k$,
denote by 
\[
s_J = s_1^{j_1} \cdot \dots \cdot s_m^{j_m} \in
\Gamma(X, L^{\otimes k})
\]
the corresponding monomial in the $s_j$, and let $D_J =
\sum j_\alpha D_\alpha$ be the divisor of $s_J$. Then
for $\ell \ge k$ multiplication by $s_J$ determines as
above a mapping
\[
\sigma_J : \II_{\ell -k} \otimes L^{\otimes -k} \lra
\II_l,
\]
and we have the following extension of Proposition
\ref{T1}:
\begin{proposition} \label{T21} \begin{enumerate}
\item[(i).] \textnormal{(cf. \cite{Lipman}.)}
If
$\ell
\ge
\min(m+k-1,n+k-1)$ then the sheaf homomorphism
\[ 
\sigma =_{\text{def}} \sum_{|J| = k} \sigma_J :
\bigoplus_{|J| = k}
  \II_{\ell -k} \otimes L^{\otimes -k} {\lra}
\II_\ell
\]
is surjective. In particular, \[\II_\ell \subset
\JJ^k.\]

\item[(ii).] Assume that $X$ is
projective, and fix  $\ell  \ge \min(m + k - 1 , n+k)$.
Let $A$ be a divisor on
$X$ such that
$A  - \ell D$ is ample (or big and nef). Then the  map
on global sections
\[ 
\bigoplus_{|J| =k} H^0(X, \OO_X(K_X + A - D_J) \otimes
\II_{\ell -k} ) \lra H^0(X, \OO_X(K_X + A) \otimes
\II_\ell) \]
induced by $\sigma$ is surjective. In particular if 
\[ s
\in H^0(X,
\OO_X(K_X + A))\]
 lies in the subspace
$
H^0(X, \OO_X(K_X + A)\otimes \II_{\ell})$, then
\[ s = \sum h_J s_J \qquad \text{for some} \ \
h_J \in H^0(X, \OO_X(K_X + A - D_J)).
\]
\end{enumerate}
\end{proposition}
\begin{proof}[Sketch of Proof]
We merely indicate the modifications required in the
proof of Proposition  \ref{T1}. Starting as before with the
surjective map of vector bundles $P \lra \OO_Y$ on $Y$,
we take $k^{\text{th}}$ symmetric powers to get $S^k
P \lra \OO_Y$. The main point is then to exhibit a
complex resolving the kernel of this map. But in fact 
there is a long exact sequence of bundles
\begin{equation}
0 \lra S^{k,1^{\times (m-1)}}P \lra \dots \lra
S^{k,1,1}P \lra S^{k,1}P \lra S^k P \lra \OO_Y \lra 0.
\tag{+}
\end{equation}
Here $S^{k, 1^{\times p}}P$ denotes the bundle formed
from $P$ via the representation of the general linear
group $GL(m, \CC)$ corresponding to the Young diagram 
$(k,1^{\times p}) = (k, 1, \dots , 1)$ ($p$ repetitions
of $1$). The existence of (+), and the fact that it
terminates where indicated, follow e.g. from
\cite{Green}, (1.a.10). Now since $P$ is a direct sum of
copies of $N^*$, it follows that
$S^{k,1^{\times p}}$ is a sum of copies of $N^{\otimes
-(k + p)}$. From this point on the argument proceeds as
before, using the twist of (+) by $Q$ in place of the
Koszul complex (*) appearing in the proof of \ref{T1}. 
\end{proof}

\larger
\vskip 5pt

\section{Distinguished Subvarieties}

In order for the results of the previous section to be
useful, one needs a criterion to guarantee that a
function $\phi$ lies in the ideal $\II_\ell$ occuring
there. We use an approach suggested by the proof of
Proposition 4.1 of \cite{Kollar1}, the idea 
being in effect to work directly on the blow-up of the
ideal sheaf $\JJ$. This naturally leads to a condition
involving  the order of vanishing of $\phi$ along
certain distinguished subvarieties of $Z$. We also give
a geometric characterization of these distinguished
subvarieties, in the spirit of \cite{VanGastel} and
\cite{FlennerVogel}, that clarifies somewhat their
connection with constructions of
\cite{Brownawell2} and \cite{Kollar2}.

We keep notation as in \S 1. Thus $X$ is a smooth
quasi-projective complex variety of dimension $n$, $\JJ
\subset \OO_X$ is an ideal sheaf defining a subscheme $B
\subset X$, and $s_j \in \Gamma(X, \JJ \otimes L)$ are
global sections generating $\JJ \otimes L$, cutting out
effective divisors $D_j$. We denote by
\[
Z = (B)_{\text{red}} = \text{Zeroes}(\sqrt{\JJ})
\]
the reduced subscheme of $X$ supported on $B$.

As above, we start by blowing up $X$ along the ideal
$\JJ$ to get
\[
\nu_0 : V_0 = \blowup_B(X) \lra X.
\]
Now let $V \lra V_0$ be the normalization of $V_0$,
with 
\[ 
\nu : V \lra X
\] the natural composition. Denote by $E$ the pull-back
to $V$ of the exceptional divisor on $V_0$, so that $E$
is an effective Cartier divisor on $V$. Then
$\JJ \cdot \OO_V = \OO_V(-E)$, and consequently 
\[
M =_{\text{def}}\  \nu^* L \  \big( -E \big) 
\]
is base-point free. Observe that since $V$ is normal,
the resolution $f : Y \lra X$ of $V_0$ introduced in \S
1 necessarily factors through a map 
\[ h : Y \lra V. \]
Moreover we have
\[
h^* E = F   \qquad, \qquad h^*M = N.
\]

Now $E$ determines a Weil divisor on $V$, say
\[
[E] = \sum_{i =1}^t r_i \ [E_i],
\]
where the $E_i$ are the irreducible components of
the support of $E$, and $r_i > 0$. Set
\[
Z_i = \nu(E_i) \subseteq X,
\]
so that $Z_i$ is a reduced and irreducible subvariety
of $X$. Remark that 
\[ 
Z_i \subset Z, \qquad \text{and} \qquad Z = \bigcup Z_i.
\]
Following Fulton and MacPherson \cite{Fulton} we call
$Z_i$ the {\it distinguished subvarieties} of $Z$, and
we refer to $r_i$ as the coefficient attached to
$Z_i$.\footnote{Strictly speaking, Fulton and
MacPherson define the distinguished subvarieties to be
the images in $X$ of the components of the exceptional
divisor of
$V_0$, but normalizing does not affect the subvarieties
that arise.} Note that two or more of the components
$E_i
\subset V$ may have the same image in $X$, i.e. there
might be coincidences among the $Z_i$, but it will be
clear that this doesn't cause any problems. (If one
wants to eliminate duplications, one could attach to
each distinct distinguished subvariety the largest
coefficient associated to it. However we prefer to
allow repetitions.) A geometric characterization of
these subvarieties is given in Proposition
\ref{gxx}

The criterion for which we are aiming is
\begin{lemma} \label{L1} Let $\II_{Z_i} \subset \OO_X$
be the ideal sheaf of $Z_i$, and denote by
$\II_{Z_i}^{<r>}$ its $r^{\text{th}}$ symbolic power,
consisting of germs of functions that have multiplicity
$\ge r$ at a general point of $Z_i$. Then for any
$\ell \ge 1$ one has the inclusion
\[
\II_{Z_1}^{<r_1 \ell>} \cap \dots \cap \II_{Z_t}^{<r_t
\ell>} \subseteq \II_\ell, \label{zz}
\]
where $\II_\ell$ is the multiplier-type ideal
introduced in \S 1.
\end{lemma}
\noi In other words, in order that a function (germ)
$\phi$ lie in $\II_\ell$, it suffices that $\phi$ have
multiplicity $\ge r_i \ell$ at a general (and hence
every) point of each of the distinguished subvarieties
$Z_i$.
\begin{corollary} \label{c22}
\begin{enumerate} 
\item[(i).] Setting $p = \min(m,n)$, one has the
inclusion
\[
\II_{Z_1}^{<r_1  p>} \cap \dots \cap \II_{Z_t}^{<r_t
 p>} \subseteq \JJ.   \]

\sbl
\item[(ii).] Assume that $X$ is projective, fix $\ell
\ge min(m , n+1)$, and let
$A$ be a divisor on $X$ such that
$A - \ell D$ is ample (or big and nef). If $s \in
\Gamma(X, \OO_X(K_X + A) )$ vanishes
to order $\ge r_i \ell$ at the general point of each
$Z_i$, then 
\[
s = \sum s_j h_j \quad \text{for some} 
\ h_j \in \Gamma(X,
\OO_X(K_X + A - D_j)).
\]
\end{enumerate}
\end{corollary}
\begin{proof} Apply 
 Proposition \ref{T1}. \end{proof}

\begin{proof}[Proof of Lemma \ref{L1}]
The assertion is local on $X$, but to avoid heavy
notation we will abusively write $X$ where we really
mean a small open subset thereof. This being said,
consider the factorization
\[ Y \overset{h}{\lra} V \overset{\nu}{\lra} X
\]
of $f : Y \lra X$, and suppose given a germ
\[ \phi \in \II_{Z_1}^{<r_1 \ell>} \cap \dots \cap
\II_{Z_t}^{<r_t
\ell>} \subset \OO_X.\]
Then $\phi$ has multiplicity $\ge r_i \ell$ at each
point of $Z_i$, and consequently $\nu^* \phi$ has
multiplicity $\ge r_i \ell$ at a general point of $E_i$
(which in particular is a smooth point of $V$). This
implies that
\[\text{ord}_{E_i}(\nu^*
\phi)
\ge r_i
\ell,\] and hence that $\text{div}(\nu^* \phi) \succeq
\ell E$. Now $F = h^* E$ and therefore
$\text{div} (f^* \phi) \succeq \ell F$. 
Since $K_{Y/X}$ is effective, this in turn implies that
\[
\text{div}(f^*\phi) + K_{Y/X} \succeq \ell F. \]
But this means exactly that 
\[
\phi \in f_* \big( \OO_Y(K_{Y/X} - \ell F) \big ) =
\II_\ell, \]
as required.  \end{proof}

\begin{example} \label{counterexample}
Here is an example to show that the factor $p =
\min(n,m)$ cannot in general be omitted from the
exponents in Corollary \ref{c22}. Fix a positive
integer $a$.  Working in
$X = \CC^2$ with coordinates $x$ and $y$, consider the
divisors defined by $s_1 = x^a$ and $s_2 = y^a$. An
explicit calculation shows that the normalized blow-up
of $X$ along the ideal $(x^a, y^a)$ is isomorphic to its
blow-up along $(x,y)$, but with exceptional divisor
$a$ times the exceptional divisor of the ``classical"
blow-up. So in this case there is a single
distinguished subvariety
$Z_1 = \{ (0,0)\}$ which appears with coefficient $r_1 =
a$. But for $a \ge 2$
\[ (x, y)^a \not \subseteq (x^a, y^a) , \quad {
i.e.} \quad   \II_{Z_1}^{<r_1>} \not \subseteq \JJ, 
\] although of course $(x,y)^{2a} \subseteq (x^a, y^a)$,
as predicted by \ref{c22}.
\end{example}

\begin{remark} \label{rtt} One can recover part (i) of
the Corollary \ref{c22} directly from the theorem of
Brian\c con-Skoda  (cf. \cite{Huneke}, Chapter 5, or
\cite{Lipman}) for regular local rings. Indeed, arguing
as above we have:
\[
  \II_{Z_1}^{<r_1
\ell>}
\cap \dots \cap \II_{Z_t}^{<r_t
\ell>}
\subset \nu_* \OO_V(-\ell E) = \overline{\JJ^\ell} . 
\]
But Brian\c con-Skoda states that
\[
\overline{\JJ^{\min(n,m)}}\subset \JJ.
\]
This suggests that in fact  the local effective
Nullstellensatz  should hold in considerably
greater algebraic generality than that which we consider
here. However it is not immediately clear in a purely
local setting how to get useful upper bounds on the
coefficients
$r_i$ in the exceptional divisor of the normalized
blow-up along the given ideal. 
\end{remark}

\vskip 5pt

\smaller

As before, the Corollary extends in a natural way to
powers of $\JJ$. In fact, Proposition \ref{T21} and the
previous Lemma yield:
\begin{corollary} \label{ppttt}
\begin{enumerate} 
\item[(i).] Setting $p = \min(m+k-1,n+k-1)$, one has the
inclusion
\[
\II_{Z_1}^{<r_1  p>} \cap \dots \cap \II_{Z_t}^{<r_t
 p>} \subseteq \JJ^k.   \]

\sbl
\item[(ii).] Assume that $X$ is projective, fix $\ell
\ge min(m +k -1, n+k)$, and let
$A$ be a divisor on $X$ such that
$A - \ell D$ is ample (or big and nef). If $s \in
\Gamma(X, \OO_X(K_X + A) )$ vanishes
to order $\ge r_i \ell$ at the general point of each
$Z_i$, then 
\[
s = \sum_{|J|=k} s_J h_J \quad \text{for some} 
\ h_J \in \Gamma \big( X,
\OO_X(K_X + A - D_J) \big). \qed
\]
\end{enumerate} 
\end{corollary}
\larger

\sbl
We conclude this section with a 
geometric characterization of the distinguished
subvarieties $Z_i \subset X$ associated to $\JJ$,
following ideas of \cite{VanGastel} and
\cite{FlennerVogel}. It shows that they are in fact
closely connected to constructions appearing in
\cite{Kollar1},
\cite{Brownawell2} and \cite{Kollar2}.  In a word, the
decomposition considered here is related to this
earlier work in much the same fashion that the
Fulton-MacPherson  intersection classes
are related to the intersection cycles constructed by
Vogel et. al. 
(which appear very explicitly in \cite{Kollar2}). 

 Let  \[ U \subset \Gamma(X , \JJ \otimes L) \]
be the $m$-dimensional subspace spanned by the
generating sections $s_1, \dots, s_m$. Given a
subspace $W \subseteq U$, set
\begin{gather*}
S_W^o = \big \{ x \in X - Z  \ \big  | \  s(x) = 0
\quad \forall s \in W \big \}, \\
S_W = \text{closure} \big (  S_W^o \big ) \subset X.
\end{gather*}
If $W \subset U$ is a general subspace of dimension
$e$, then $S_W$ is an algebraic subset of $X$ of pure
dimension $n - e$. 
\begin{proposition} \label{gxx}
Let $T \subset X$ be an irreducible subvariety of
dimension $d \le n-2$, and consider the
$(d+1)-dimensional$ subsets $S_W \subset X$ for   
$W \subset U$ a general subspace of dimension $n - d -
1$. Then $T$ is distinguished if and only if
\[ T \subset S_W \]
for all sufficiently general $W$.
\end{proposition}
\noi For example consider the case $d = 0$, so that $T$
is a single point. Then the subsets $S_W$ appearing in
the Proposition are curves, and the assertion is
that the distinguished points are exactly the common
intersection points of this family of curves. 

\begin{proof}[Proof of Proposition \ref{gxx}]

 The vector space $U \subset \Gamma(X, \JJ
\otimes L)$ is isomorphic in the natural way to a
subspace
$U^\prime \subset \Gamma(V, M)$ generating
$M = \nu^*L \ (-E)$, and in the sequel we identify $U$
and $U^\prime$. We will consider the maps
\begin{equation}
\begin{CD}
V @>\phi>> \PP^{m-1} \\
@V\nu VV \\
X
\end{CD}, \tag{+}
\end{equation}
$\phi$
being the morphism defined by $U^\prime$ (or $U$), so 
that $\phi^* \OO_{\PP^{m-1}}(1) = M$. Thus
$\PP^{m-1}$ is  the projective space of
one-dimensional quotients of $U$, and a subspace $W
\subset U$ of dimension $e$ corresponds to a linear
subspace $L_W \subset \PP^{m-1}$ of codimension $e$.

We claim first  that all of the fibres of $\nu$
map finitely to $\PP^{m-1}$. In fact, 
the blow-up $V_0 = Bl_B(X)$ is the closure
of the graph of the rational map $X \dashrightarrow
\PP^{m-1}$ determined by the $s_j \in \Gamma(X,L)$ (cf.
\cite{Fulton}. Chapter 4, \S 4).
Thus
$V_0$ sits naturally as a subvariety
\[
V_0 \subset X \times \PP^{m-1},
\]
and in particular  its normalization
$V$ maps finitely to $X \times \PP^{m-1}$ via the
morphism determined by (+). Hence the fibres of $V$ over
$X$ are indeed finite over
$\PP^{m-1}$, as claimed. It follows in particular that
if $E_i$ is a component of the exceptional divisor $E =
\nu^{-1}(Z)$ in $V$, and if
$Z_i = \nu(E_i) \subset X$ is the corresponding
distinguished subvariety of $X$, then all of the fibres
of $E_i \lra Z_i$ map finitely to $\PP^{m-1}$.

We claim next that if
$W \subset U$ is a sufficiently general subspace of
dimension $1
\le n - d - 1 \le n - 1$, then 
\begin{equation} S_W = \nu \big( \phi^{-1}(L_W) \big ).
\tag{*}
\end{equation}
In fact, the two sides of (*) evidently agree away from
$Z$. So to verify that they actually coincide, it
suffices to show that no irreducible component of
$\phi^{-1}(L_W)$ is contained in the exceptional
divisor $E = \nu^{-1}(Z)$. To this end, let $E_i$
denote an irreducible component of $E$, so
that  $E_i$ has dimension $n - 1$. Then for 
sufficiently general $W$, either $E_i
\cap \phi^{-1}(L_W) = \emptyset$ or else
\[ 
 \dim \big ( E_i \cap \phi^{-1}(L_W) \big ) = (n-1) -
(n-d-1) = d.\] On the other hand, $\phi^{-1}(L_W)$
itself is either empty or of pure dimension $d + 1$,
and so indeed no component of $\phi^{-1}(L_W)$ is
contained in the support of $E$.

 Now fix an irreducible subvariety $T \subset X$ of
dimension $d \le n - 2$.  Then $T$ is distinguished iff
$\nu^{-1}(T)$ contains at least one irreducible
component of dimension $n - 1$ which dominates $T$.
Setting $F_t = \nu^{-1}(t)$, this is in turn equivalent
to the condition that for general
$t \in T$:
\begin{equation} \dim F_t \ge n - d - 1. \notag
\end{equation}
We have noted already that $\phi$ restricts to a finite
mapping on each of the fibres $F_t$, and since $\nu$ is
proper $F_t$ is complete. Therefore $\phi(F_t) \subset
\PP^{m-1}$ is a Zariski-closed subset having the same
dimension as $F_t$. Thus
$\dim F_t \ge n - d - 1$ if and only if $\phi(F_t)$
meets any linear space $L \subset \PP^{m-1}$ of
codimension $n - 1 - d$, i.e. iff
\begin{equation}
F_t \cap \phi^{-1}(L_W) \ne \emptyset \tag{**}
\end{equation}
for general $W \subset U$ of dimension $n - d-1$. But
(**) holds for general $t \in T$ iff 
\[ \nu \big( \phi^{-1}(L_W) \big ) \supseteq T. \]
The Proposition then follows from (*).
\end{proof}

\begin{remark} It would be interesting to have a
geometric characterization of the coefficients $r_i$
attached to the distinguished subvarieties $Z_i$.
\end{remark}

\vskip 5pt

\section{Degree Bounds}

The only remaining point is to prove a Brownawell-type
bound on the degrees of the distinguished
subvarieties $Z_i$. There are general positivity
theorems for Fulton-MacPherson intersection classes, 
as developed e.g. in \cite{FL}, 
lurking here. However it is easiest to bypass these
results in the case at hand. 

We keep notation as in the previous sections.  Thus $X$ 
is a smooth complex variety of dimension $n$, $\JJ
\subset \OO_X$ is an ideal sheaf defining a subscheme
$B \subset X$, $Z = B_{\text{red}}$ is the
corresponding reduced algebraic subset of $X$, and
$s_1, \dots, s_m \in \Gamma(X, \JJ \otimes L)$ are  
global sections generating $\JJ
\otimes L$. We continue to denote by
$Z_1,
\dots, Z_t
\subset X$  the distinguished subvarieties determined
by $\JJ$, and by
$r_i > 0$ the coefficients attached to them. 

\begin{proposition}  \label{pp3}
Assume that $X$ is projective and that $L$ is
nef. Then
\begin{equation} 
\sum_{i = 1}^t  \ r_i  \cdot \deg_L(Z_i) \ \le \
\deg_L(X) \ = \ \int_X c_1(L)^n.
\tag{*}
\end{equation}
\end{proposition}
\newcommand{ \Lt} {\tilde L}
\begin{proof}
The given sections $s_j \in \Gamma(X, \JJ \otimes L)$
determine in the natural way sections $s_j^\prime \in
\Gamma(V, M)$ generating $M = \nu^*L (-E)$. We consider
as in the proof of Proposition
\ref{gxx} the corresponding morphism 
\[ \phi : V \lra
\PP^{m-1},  \]
 so  that $\phi^* \OO_{\PP^{m-1}}(1)
= M$. Recall from that proof that $\phi$ is finite on
all the fibres of $\nu$ (the point being that $V_0 =
\text{Bl}_B(X)$ embeds as a subvariety of $X \times
\PP^{m-1}$, and hence that $V$ maps finitely to this
product). In particular, if $E_i$ is a component of the
exceptional divisor $E$ in $V$, and if $Z_i = \nu(E_i)$
is the corresponding distinguished subvariety, then the
restriction of $M$ to any of the fibres of $E_i \lra
Z_i$ is ample.

Now denote by  $\Lt = \nu^* L$ the
pull-back of
$L$ to
$V$. Noting that $\int_V c_1(M)^n \ge 0$, and recalling
that
$[E] = c_1(\Lt) - c_1(M)$, we have:
\[
\begin{aligned}
 \deg_L(X) &= \int_V c_1(\Lt)^n \\
 &\ge \int_V \big( c_1(\Lt)^n - c_1(M)^n \big ) \\
 &= \int_V \Big( c_1(\Lt) - c_1(M) \Big) \Big(  \sum_{j
= 0}^{n-1} c_1(\Lt)^{j} c_1(M)^{n-1-j} \Big) \\
&= \int_{[E]} \Big(  \sum_{j
= 0}^{n-1} c_1(\Lt)^{j} c_1(M)^{n -1-j} \Big)  \\ 
&= \sum_{i=1}^t r_i \cdot \int_{E_i}\Big(  \sum_{j
= 0}^{n-1} c_1(\Lt)^{j} c_1(M)^{n-1-j} \Big) \\
&\ge \sum_{i = 1}^t r_i \cdot \int_{E_i}
c_1(\Lt)^{\dim(Z_i)}c_1(M)^{n-1-\dim(Z_i)}, \\
\end{aligned}
\]
where in the last step we have used  that 
\[ 
\int_{E_i} c_1(\Lt)^j c_1(M)^{n-1-j} \ge 0 \ \ 
\text{for all $j$}
\]
 thanks to the fact that $\Lt$ and $M$ are nef. Now the
restriction to $E_i$ of $c_1(\Lt)^{\dim(Z_i)}$ is
represented (say in rational cohomology) by 
$\deg_L(Z_i)$ general fibres of the map
$E_i \lra Z_i$. Moreover as we have noted the
restriction of $M$ to each of these fibres is ample, 
and hence each has positive $M$-degree. Therefore
\[
\int_{E_i}
c_1(\Lt)^{\dim(Z_i)}c_1(M)^{n-1-\dim(Z_i)} \ge
\deg_L(Z_i),
\]
and the Proposition follows.
\end{proof}

\vskip 5pt

The proof of the Theorem stated in the Introduction is
now complete. The degree bound just established
combined with Corollary \ref{ppttt} also give an
analogous statement involving higher powers of $\JJ$.

\begin{remark} One can obtain a slight strengthening of
Proposition \ref{pp3} by taking into account a further
geometric invariant. Specifically, denote by $\mu$
the number of intersection points away from $Z$ of
$n$ general divisors in the linear series spanned by
the $D_j$ (and set $\mu = 0$ if $m \le n$).
Equivalently, with notation as at the end of
\S 2, $\mu = \# S_W$ where $W \subset U$ is a general
subspace of dimension $n$. Then in the situation of
\ref{pp3} the calculations just completed show that in
fact:
\[
\sum_{i = 1}^t  \ r_i  \cdot \deg_L(Z_i) \ \le \
\deg_L(X) - \mu.
\]
Indeed, simply observe that with notation as in the
previous proof:
\[ \mu = \int_V c_1(M)^n .\]

\end{remark}

\begin{remark} In the statement of the Theorem
appearing in the Introduction, we assumed for
simplicity that the line bundle $L$ is ample. In fact,
the only positivity  used in the proof is the nefness
of $L$, which comes into Proposition \ref{pp3}.
However to get a non-trivial assertion, one wants to
avoid the possibility that the $L$-degrees appearing
there might be zero. Perhaps then the most  natural 
hypothesis for  the Theorem is that $L$ is nef, and
that its restriction to the zero-locus $Z$ is ample. By
the same token, in statement (iii) of the main Theorem,
it is sufficient to suppose that $A - (n+1)D$ is big
and nef. 
\end{remark}

\begin{remark} In the work of Koll\'ar  \cite{Kollar1}
and others on projective space, one allows the degrees
of the defining equations to differ. One can generalize
the results here to the case where the divisors $D_j$
lie in different linear series by imposing the condition
that $\OO_X(D_j - D_k)$ be base-point-free for $j \le
k$. However this makes the arguments a little more
technical and less transparent, and we do not address
this extension here. 
\end{remark}

\end{document}